\documentclass[review]{elsarticle}

\usepackage{lineno,hyperref}
\modulolinenumbers[5]
\usepackage{amsmath,amssymb,fixmath}
\usepackage{epsfig,epstopdf}
\usepackage{booktabs}
\usepackage{multirow}
\usepackage{subfig}
\usepackage{float}
\journal{13th IEEE Colloquium on Signal Processing \& its Applications}
\begin{document}
	
	\begin{frontmatter}		
		\title{A Robust Variable Step Size Fractional Least Mean Square (RVSS-FLMS) Algorithm}
		
		\author[mymainaddress,b]{Shujaat Khan}
		\ead{shujaat@kaist.ac.kr}
		
		\author[b]{Muhammad Usman}
		\ead{musman@iqra.edu.pk}
		
		\author[c,d]{Imran Naseem\corref{mycorrespondingauthor}}
		\cortext[mycorrespondingauthor]{Corresponding author}
		\ead{imran.naseem@uwa.edu.au}
		
		\author[d]{Roberto Togneri}
		\ead{roberto.togneri@uwa.edu.au}
		
		\author[e]{Mohammed Bennamoun}
		\ead{mohammed.bennamoun@uwa.edu.au}
		
		\address[mymainaddress]{Department of Bio and Brain Engineering, Korea Advanced Institute of Science and Technology (KAIST), Daejeon, Republic of Korea.}
		\address[b]{Faculty of Engineering Science and Technology (FEST), Iqra University, Defence View, Karachi-75500, Pakistan.}
		\address[c]{College of Engineering, Karachi Institute of Economics and Technology,	Korangi Creek, Karachi 75190, Pakistan.}
		\address[d]{School of Electrical, Electronic and Computer Engineering, The University of Western Australia, 35 Stirling Highway, Crawley, Western Australia 6009, Australia.}
		\address[e]{School of Computer Science and Software Engineering, The University of Western Australia, 35 Stirling Highway, Crawley, Western Australia 6009, Australia.}

\begin{abstract}
In this paper, we propose an adaptive framework for the variable step size of the fractional least mean square (FLMS) algorithm.  The proposed algorithm  named the robust variable step size-FLMS (RVSS-FLMS), dynamically updates the step size of the FLMS to achieve high convergence rate with low steady state error.  For the evaluation purpose, the problem of system identification is considered.  The experiments clearly show that the proposed approach achieves better convergence rate compared to the FLMS and adaptive step-size modified FLMS (AMFLMS).
\end{abstract}
\end{frontmatter}
\section{Introduction}
The standard integral and derivative are essential tools of calculus used by the professionals dealing with the natural or artificial systems. The fractional order calculus (FOC) however extends the traditional definitions which are limited to integer values.  The FOC is as old as classical calculus, however, the applications of the FOC in the areas of engineering and research have become more popular in the last decade. New areas are emerging as both theoretical and practical aspects of the FOC operators are now well established. Fractional calculus has turned out to be an important tool in various fields \cite{FCA5}, \cite{FCA6}. 

Liouville, Reimann and Leibniz, amongst others, are considered to be the pioneers of fractional calculus \cite{liou}.  In the late $19^{th}$ century Heaviside made the use of fractional derivatives in transmission lines and gave the solution for the diffusion equation.  In the recent years fractional calculus has been applied to the problem of path tracking in autonomous vehicles \cite{FCA1}. In \cite{FCA2} the concept of pulses propagating through porous media is proposed with an experimental proof for a theoretical model based on fractional derivatives. Fractional calculus is also applied in the theory of viscoelasticity \cite{FCA3}.  In the field of image processing fractional differentiation has been used for a number of tasks such as edge detection \cite{FCA4}.  In the domain of signal processing FOC has found a number of useful applications. Adaptive filtering for instance, has been the main focus of research in this context.  In \cite{LMS12} an FOC based least mean square algorithm named the fractional least mean square (FLMS) is proposed.  The FLMS is compared with the least mean square (LMS) algorithm for the problem of system identification.

In the recent past, various modifications have been proposed for the conventional FLMS \cite{LMS14, LMS16, LMS17,FNSAF,FONFLMS,NFVFLMS,MCCFLMS,RVPFLMS}.  In this paper a modified FLMS based algorithm is proposed.  The proposed method is named the robust variable step size fractional least mean square (RVSS-FLMS) and is inspired by the variable learning rate LMS (RVSS-LMS) \cite{RVSS1}.  The RVSS-FLMS algorithm dynamically adapts the learning rate of FLMS using the instantaneous error energy to achieve high convergence rate and low steady state error.  The proposed algorithm is evaluated on the problem of system identification and is compared with the conventional FLMS and adaptive step-size modified fractional least mean square (AMFLMS) algorithm.  The rest of the paper is organized as follows : In section \ref{FLMS} a brief literature review of modified least mean square algorithms is provided.  In section \ref{RVSSFLMS} the proposed RVSS-FLMS algorithm is discussed followed by the experiments in section \ref{Experimental}.  The paper is concluded in section \ref{Conclusion}.

\section{Least Mean Square}\label{FLMS}
The statistics in the LMS are estimated continuously, therefore it is an adaptive filter belonging to the group of stochastic gradient methods. Extensive research is done towards optimization of the LMS algorithm by numerous researchers \cite{LMS14, LMS13, LMS11, LMS9, LMS15}.  LMS is applied in diversified applications such as plant identification \cite{LMS12}, noise cancellation \cite{NCLMS}, echo cancellation \cite{ECLMS}, ECG signal analysis \cite{ECGLMS}, time series prediction \cite{TSPLMS} etc.  Achieving fast convergence and low steady state error is a challenge, various researchers have proposed different solutions for this.  One of the disadvantages of the LMS is that it is sensitive to the scaling of its input \cite{NLMS}.  A variant of the LMS algorithm the normalised least mean squares filter (NLMS) solves this problem through normalization.  In \cite{LMS1}, an optimized  NLMS filter is proposed to achieve a tradeoff between the convergence rate and tracking.  In \cite{LMS2} q-LMS is proposed which utilizes the q-gradient from the Jackson's derivative so that secant of the cost function is computed instead of the tangent.  The algorithm therefore takes larger steps towards the optimum solution and therefore achieves the higher convergence rates.  For the optimal processing of complex signals, linear modelling in complex numbers is utilized by \cite{LMS3}. For the complex domain adaptive filtering an augmented complex least mean square (ACLMS) is proposed. 

In \cite{LMS15} fractional least mean square (FLMS) is proposed for parameter estimation of input non-linear control autoregressive (INCAR) or Hammerstein non-linear controlled auto regression models based on fractional signal processing approach. Unknown parameters are estimated by applying the FLMS algorithm having different step sizes.  The algorithm provides better convergence results when compared to the volterra least mean square (VLMS) and kernel least mean square (KLMS). A modified structure of the FLMS, for the prediction of chaotic and non-stationary time series, named the modified fractional least mean square (MFLMS) is proposed in \cite{LMS16}.  The algorithm incorporates the adjustable gain parameter thereby avoiding some complex calculations and with reduced computational expense.  The MFLMS is tested on stationary and non-stationary time series with different noise levels and has shown improved results compared to the LMS and FLMS. In \cite{LMS17} the adaptive weight gain parameters are incorporated by implementing a gradient-based approach on the variable learning scheme.  Implementation of this approach changes the nature of the order of fractional derivative in MFLMS from fixed to adaptive.  The adaptive nature of the adjustable parameters may help the filters in handling the problems with non-linear nature.  They proposed adaptive step-size modified fractional least mean square (AMFLMS) and demonstrated the improvements by comparing their algorithm with the LMS, FLMS and MFLMS. 

In this research we propose to utilize the concept of robust variable step size (RVSS) \cite{RVSS1} for the variable learning rate of the FLMS algorithm.  The proposed scheme is robust and computationally less expensive.  For the performance evaluation we consider the problem of the system identification and compare our results with the FLMS and AMFLMS algorithm.

\section{Proposed RVSS FLMS}\label{RVSSFLMS}
In \cite{LMS12} the weight update equation for the fractional LMS (FLMS) is given as:
\begin{eqnarray}\label{weightupdate1}
	w_{k}(n+1)= w_{k}(n)- \nu \dfrac{\partial J(n)}{\partial w_{k}} - \nu_{f}\left(\dfrac{\partial}{\partial w_{k}}\right)^f J(n)
\end{eqnarray}

where $w_{k}(n)$ is the weight of the $k^{th}$ tap at the $n^{th}$ iteration, $f$ is the fractional power of derivative, and $\nu$ and $\nu_{f}$ are the step sizes.

\begin{eqnarray}\label{cost}
	J(n)= \dfrac{1}{2}e(n)^2 = \dfrac{1}{2}(d(n)-y(n))^{2}
\end{eqnarray}

$J(n)$ is the cost function defined in Eq. (\ref{cost}), where $e(n)$ is the instantaneous error between the desired output $d(n)$ and the estimated output $y(n)$ at the $n^{th}$ iteration.

Now, the  $\dfrac{\partial J(n)}{\partial w_{k}}$ is defined as :

\begin{eqnarray}\label{costd}
	\dfrac{\partial J(n)}{\partial w_{k}} = \dfrac{\partial J(n)}{\partial e} \dfrac{\partial e(n)}{\partial y} \dfrac{\partial y(n)}{\partial w_{k}}
\end{eqnarray}

For the fractional derivative term $\left(\dfrac{\partial}{\partial w_{k}}\right)^f J(n)$, the methods of chain rule can be defined in various ways. The modified Rieman-Lioville based definition of the fractional chain rule i.e $ D_x^\nu f(g(x)) = (D_g^1f(g))_{g=g(x)} D_x^\nu g(x) $, as suggested in \cite{FOCCR1}, is one of the well known techniques.  This definition has shown good results in various applications of fractional calculus. For example, in \cite{FOCCR1} for the fractional Taylor series of non differentiable functions, \cite{FOCCR5} for the derivation of further fractional order derivatives, \cite{FOCCR7} for system modeling and others  \cite{FOCCR2,FOCCR3,FOCCR4,FOCCR6}.

Using the modified Rieman-Lioville chain rule for fractional derivatives, the $\left(\dfrac{\partial}{\partial w_{k}}\right)^f J(n)$ term is defined as :

\begin{eqnarray}\label{costdf}
	\left(\dfrac{\partial}{\partial w_{k}}\right)^f J(n) = \dfrac{\partial J(n)}{\partial e} \dfrac{\partial e(n)}{\partial y} \left(\dfrac{\partial}{\partial w_{k}}\right)^f y(n)
\end{eqnarray}

solving Eq. \eqref{costd} results in :
\begin{eqnarray}\label{costd1}
	\dfrac{\partial J(n)}{\partial w_{k}} = -e(n)x(n)
\end{eqnarray}

According to the Rieman-Lioville fractional derivative method :
\begin{eqnarray}
	D^v f(t)= \dfrac{1}{\Gamma(n-v)}\left(\dfrac{d}{dt}\right)^n \int\limits_{0}^{t}(t-\tau)^{(n-v-1)} f(\tau)d\tau
\end{eqnarray}
\begin{eqnarray}\label{RemLov}
	D^v(t-a)^\alpha = \dfrac{\Gamma(1+\alpha)}{\Gamma(1+\alpha-v)}(t-a)^{(\alpha-v)}
\end{eqnarray}
where $\alpha-v+1 > 0 , D$ is the differential operator and $v$ is a real number defining the fractional power of the derivative.

Using Eq. \eqref{RemLov}, Eq. \eqref{costdf} can be reduced to:

\begin{eqnarray}\label{costdf1}
	\left(\dfrac{\partial}{\partial w_{k}}\right)^f J(n) = -e(n)x(n)\dfrac{w_{k}^{(1-f)}(n)}{\Gamma(2-f)}
\end{eqnarray}	 

Using Eq. \eqref{costd1} and  \eqref{costdf1} the weight update Eq. \eqref{weightupdate1} becomes : 

\begin{equation}\label{weightupdate2}
	w_{k}(n+1)=w_{k}(n) + \nu e(n)x(n)  + \nu_{f}e(n)x(n)\dfrac{w_{k}^{1-f}(n)}{\Gamma(2-f)}
\end{equation}
%
%
Note that the result in Eq. (\ref{weightupdate2}) is similar to the FLMS algorithm results \cite{LMS16, LMS17, LMS15}, obtained from the alternative method of fractional derivatives to avoid the chain rule.

For simplicity, in Eq. \eqref{weightupdate2}, we consider $\nu_{f} = \nu \Gamma(2-f)$ and this results in :

\begin{equation}
	w_{k}(n+1)=w_{k}(n) + \nu e(n)x(n)  + \nu e(n)x(n) w_{k}^{1-f}(n)
\end{equation}

\begin{equation}
	w_{k}(n+1)=w_{k}(n) + \nu e(n)x(n) \left(1 + w_{k}^{1-f}(n)\right)
\end{equation}

For the time varying step size $\nu$ can be replaced by $\nu(n)$

\begin{equation}
	w_{k}(n+1)=w_{k}(n) + \nu(n) e(n)x(n) \left(1 + w_{k}^{1-f}(n)\right)
\end{equation}

For the adaptation of learning rate, we proposed to use the concept of error energy correlation of RVSS-LMS algorithm \cite{RVSS1}.  The update rule for the time varying learning rate $\nu(n)$ in the proposed robust variable step size-FLMS (RVSS-FLMS) is defined as:	
\begin{eqnarray}
	\nu (n+1)= \beta \nu(n)+ \gamma p^2(n)
\end{eqnarray}	

where $(0 < \beta < 1)$, $(\gamma > 0)$, $p(n)$ is the average error energy correlation and $\nu(n+1)$ is set to $\nu_{min}$ or $\nu_{max}$

\begin{equation}
	p(n) = \alpha p(n - 1) + (1 - \alpha) e(n) e(n - 1)
\end{equation}

The positive constant $ \alpha \; (0 < \alpha < 1) $ is a weighting parameter that governs the averaging time constant. It is referred to as the forgetting factor.  The limits on $ \nu (n + 1) $ are given by :

\begin{equation}
	\nu (n + 1) = \left\{ \begin{array}{rcl}
		\nu_{\max} & \mbox{if} & \nu (n + 1) > \nu_{\max} \\ 
		\nu_{\min} & \mbox{if} & \nu (n + 1) < \nu_{\min} \\ 
		\nu(n + 1) & & otherwise
	\end{array}\right.
\end{equation}

where $ \nu_{\max} > \nu_{\min} > 0 $. 

\section{Simulation Setup and Results}\label{Experimental}

For the evaluation of the proposed method the problem of plant identification is considered.  The modeling of a linear system is an interesting problem.  Sometimes it is desired to model the system as a linear filter to avoid unnecessary complexity. The adaptive learning methods like LMS, achieves good performance in this regard \cite{LMS12,LMS13,LMS15,LMS1}.  To evaluate the efficacy of the proposed robust variable step size-FLMS (RVSS-FLMS), we consider a linear system, shown in Fig. \ref{plant1}:
\begin{figure}[h!]
	\begin{center}
		\centering 
		\includegraphics*[scale=0.40,bb=-10 -10 600 290]{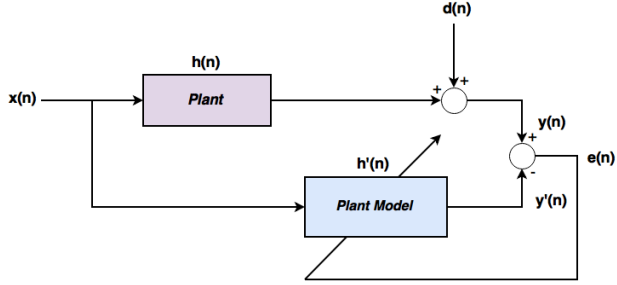}
	\end{center}
	\caption{System identification using adaptive learning algorithm.}
	\label{plant1}
\end{figure}
\begin{equation}\label{planteq}
	y(n)=a_1x(n)+a_2x(n-1)+a_3x(n-2) + d(n)
\end{equation}  

Eq. (\ref{planteq}) shows the mathematical model of the system, where $x(n)$ is the input and $y(n)$ is the output of the system, $d(n)$ is the disturbance model assumed to be $\mathcal{N}(0,\sigma_d^2)$, $a_i$'s represent the polynomial coefficients depicting the zeros of the system.  For the said experiment, $x(n)$ is a binary shift keying (BPSK) modulated random sequence comprising $600$ samples.  In Fig. \ref{plant1}, the impulse response of the system is defined by $h(n)$ while  $\hat{y}(n)$ is the estimated output, $\hat{h}(n)$ is the estimated impulse response and $e(n)$ is the estimation error.  The simulation parameters selected for the experiments are:  $a_1=0.9$, $a_2=0.3$, $a_3=-0.1$.  The experiments are performed on four noise levels with the following SNR values $10$ dB, $20$ dB, $30$ dB, and $40$ dB.

For the model design, the number of taps were chosen to be $3$.  The weights were initialized with a small value of $1\times10^{-20}$. For the fractional LMS (FLMS), adaptive step-size modified fractional least mean square (AMFLMS) and RVSS-FLMS the initial values of  $\nu(n)$ and $\nu_{f}(n)$ for $n=0$ set to $1\times10^{-4}$ and the fractional power $f(n)$ for $n=0$ was taken to be $0.5$.  For the AMFLMS, $\Psi_{k}$, $\Psi_{k,f}$ and $\Phi_{k,f}$ were initialized with a small value of $1\times10^{-15}$.  The mixing parameter $\beta$ of AMFLMS was set to $0.5$ and the adaptation rates $\alpha$, $\alpha_f$ and $\alpha_{fr}$ were all set to $1\times10^{-13}$.  For the proposed RVSS-FLMS the values of $\alpha$, $\beta$ and $\gamma$ were all set to $0.5$. The value of $\nu_{max}$ and $\nu_{min}$ are chosen to be $3\times10^{-4}$ and $1\times10^{-4}$ respectively.

\begin{figure}[h!]
	\begin{center}
		\centering 
		\includegraphics*[scale=0.52,bb=50 185 560 610]{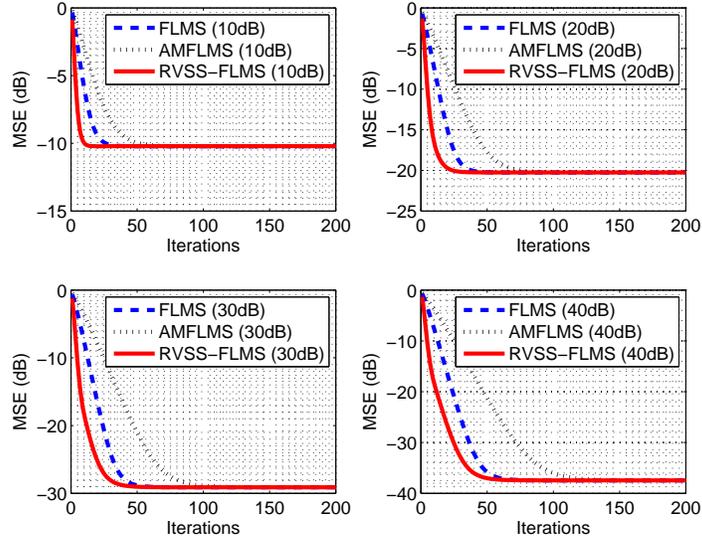} 
	\end{center}
	\caption{System identification:  The MSE curves for various methods.}
	\label{planterror}
\end{figure}

The experiments were performed on four noise levels and the mean square error (MSE) curves are depicted in Fig. \ref{planterror}.  For the SNR values of $10$ dB, $20$ dB, $30$ dB and $40$ dB the RVSS-FLMS produced the best performance, achieving an MSE of $-10.22$ dB, $-20.26$ dB, $-29.70$ dB and $-37.68$ dB in around 20, 38, 60 and 65 iterations respectively. The conventional FLMS converged to the MSE value of $-10.21$ dB, $-20.25$ dB, $-29.71$ dB and $-37.58$ dB in around 45, 60, 80 and 90 iterations respectively.  Whereas the AMFLMS converged to the MSE value of $-10.22$ dB, $-20.27$ dB, $-29.71$ dB and $-37.60$ dB in around 80, 100, 120 and 140 iterations respectively.  The convergence rate of the proposed RVSS-FLMS algorithm is therefore found to be superior the FLMS and the AMFLMS algorithms.

\begin{figure}[h!]
	\begin{center}
		\centering 
		\includegraphics*[scale=0.52,bb=50 185 560 610]{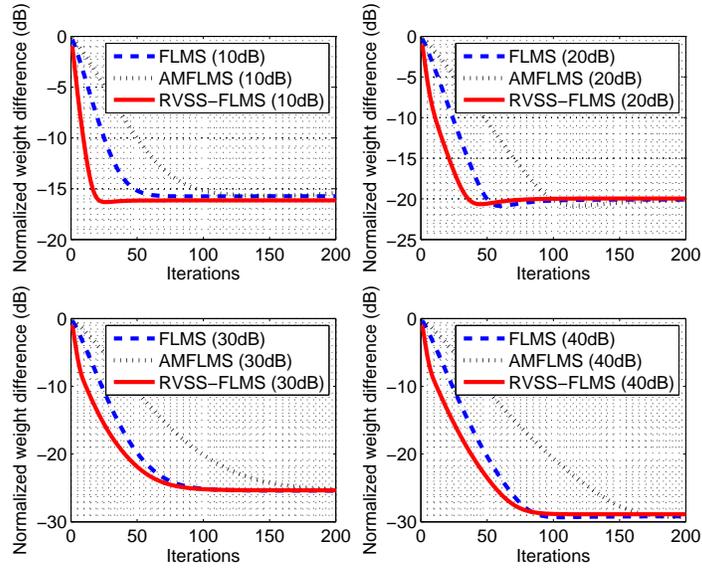} 
	\end{center}
	\caption{System identification:  The normalized weight difference curves for various methods.}
	\label{plantNWD}
\end{figure}

For the fair evaluation of the proposed algorithm, another well known performance measurement parameter called the normalized weight difference (NWD), was used.  The NWD is a comparative measure between the actual impulse response and the identified weights. The experiments were performed on four noise levels for which the NWD curves are depicted in Fig. \ref{plantNWD}.  For the SNR value of $10$ dB, $20$ dB, $30$ dB and $40$ dB the RVSS-FLMS produced high convergence and low steady state error, achieving an NWD value of $-16.21$ dB, $-20.00$ dB, $-25.06$ dB and $-28.52$ dB in approximately 35, 60, 95 and 100 iterations respectively.  The FLMS converges to the NWD value of $-15.81$ dB, $-20.09$ dB, $-25.06$ dB and $-28.91$ dB in approximately 60, 85, 120 and 105 iterations respectively.  Whereas the AMFLMS converges to the NWD value of $-15.67$ dB, $-20.12$ dB, $-25.00$ dB and $-29.02$ dB in approximately 120, 175, 200 and 200 iterations respectively.  The proposed algorithm has therefore shown to comprehensively outperform the FLMS and AMFLMS algorithms.

For the purpose of time-complexity comparison of the proposed algorithm with the FLMS and AFLMS, we investigated the training time of FLMS, AMFLMS and RVSS-FLMS for 200 iterations.  The proposed method utilized 2.05 seconds whereas the FLMS and AMFLMS took 1.87 seconds and 11.06 seconds respectively.  All the experiments were conducted using Matlab on an Intel(R) Core(TM) i5-4690
CPU @ 3.5GHz machine with 16GB memory.  Although the proposed algorithm took a few milliseconds more than the FLMS but the experiments clearly show that the proposed robust approach dynamically adapts the learning rate to achieve the minimum steady state error in lesser number of iterations.  The summary results for all simulations is showed in Table \ref{summary_table}.

\begin{table*}[]
	\caption{Results for plant identification problem}
	\label{summary_table}
	\resizebox{\textwidth}{!}{%
		\begin{tabular}{l|cccc|cccc|c}
			\multirow{1}{*}{\textbf{Algorithm}} & \multicolumn{4}{|c|}{\textbf{(Iterations, MSE) for SNR values}}           & \multicolumn{4}{|c|}{\textbf{(Iterations, NWD) for SNR values}}             & \multirow{1}{*}{\textbf{Time (seconds)}} \\
			& \textit{10 dB}  & \textit{20 dB}   & \textit{30 dB}   & \textit{40 dB}  & \textit{10 dB}   & \textit{20 dB}   & \textit{30 dB}   & \textit{40 dB}   &                                      \\
			\midrule
			\textit{\textbf{FLMS}}              & (45, -10.21 dB) & (60, -20.25 dB)  & (80, -29.71 dB)  & (90, -37.58 dB) & (60, -15.81 dB)  & (85, -20.09 dB)  & (120, -25.06 dB) & (105, -28.91 dB) & 1.87                                     \\
			\textit{\textbf{AMFLMS}}            & (80, -10.22 dB) & (100, -20.27 dB) & (120, -29.71 dB) & (140, -37.6 dB) & (120, -15.67 dB) & (175, -20.12 dB) & (200, -25.00 dB) & (200, -29.02 dB) & 11.06                                    \\
			\textit{\textbf{RVSS-FLMS}}         & (20, -10.22 dB) & (38, -20.26 dB)  & (60, -29.70 dB)  & (65, -37.68 dB) & (35, -16.21 dB)  & (60, -20.00 dB)  & (95, -25.06 dB)  & (100, -28.52 dB) & 2.05                              
	\end{tabular}}
\end{table*}

\section{Conclusion}\label{Conclusion}
In this research an adaptive least mean square algorithm based on fractional derivative is proposed. In particular, the step size of the fractional LMS (FLMS) is made adaptive using the concept of robust variable step size.  For this purpose the concept of robust variable step size is utilized.  The performance evaluation of the proposed approach is done by implementing it for the problem of system identification with different noise levels.  The results of the proposed algorithm are compared to the FLMS and adaptive step-size modified fractional least mean square (AMFLMS) approaches.  The proposed algorithm attains better convergence rate and steady-state error and is therefore found to be superior to the FLMS and AMFLMS algorithms.

	

\section*{References}
\bibliography{ref1}

\end{document}